\pdfoutput=1
\RequirePackage{ifpdf}
\ifpdf 
\documentclass[pdftex]{sigma}
\else
\documentclass{sigma}
\fi

\numberwithin{equation}{section}
\newtheorem{Theorem}{Theorem}[section]
 { \theoremstyle{definition}
\newtheorem{Definition}[Theorem]{Definition}}

\newcommand{\R}{\mathbb{R}}
\newcommand{\Z}{\mathbb{Z}}
\newcommand{\h}{\mathbb{H}}
\newcommand{\cL}{\mathcal{L}}
\def \Gal{\operatorname{Gal}}
\def \DGal{\operatorname{DGal}}
\def \Aut{\operatorname{Aut}}
\def \GL{\operatorname{GL}}
\def \SO{\operatorname{SO}}
\def \SL{\operatorname{SL}}
\def \SU{\operatorname{SU}}
\def \U{\operatorname{U}}
\def \Sp{\operatorname{Sp}}

\begin{document}
\allowdisplaybreaks

\newcommand{\arXivNumber}{1403.3226}

\renewcommand{\PaperNumber}{100}

\FirstPageHeading

\ShortArticleName{Picard--Vessiot Extensions of Real Differential Fields}

\ArticleName{Picard--Vessiot Extensions of Real Differential Fields}

\Author{Teresa CRESPO~$^\dag$ and Zbigniew HAJTO~$^\ddag$}

\AuthorNameForHeading{T.~Crespo and Z.~Hajto}

\Address{$^\dag$~Departament de Matem\`atiques i Inform\`atica, Universitat de Barcelona,\\
\hphantom{$^\dag$}~Gran Via de les Corts Catalanes 585, 08007 Barcelona, Spain}
\EmailD{\href{mailto:teresa.crespo@ub.edu}{teresa.crespo@ub.edu}}

\Address{$^\ddag$~Faculty of Mathematics and Computer Science, Jagiellonian University,\\
\hphantom{$^\ddag$}~ul.~Prof. S. {\L}ojasiewicza~6, 30-348 Krak\'ow, Poland}
\EmailD{\href{mailto:zbigniew.hajto@uj.edu.pl}{zbigniew.hajto@uj.edu.pl}}

\ArticleDates{Received July 04, 2019, in final form December 22, 2019; Published online December 24, 2019}

\Abstract{For a linear differential equation defined over a formally real differential field~$K$ with real closed field of constants~$k$, Crespo, Hajto and van der Put proved that there exists a~unique formally real Picard--Vessiot extension up to $K$-differential automorphism. However such an equation may have Picard--Vessiot extensions which are not formally real fields. The differential Galois group of a Picard--Vessiot extension for this equation has the structure of a linear algebraic group defined over $k$ and is a $k$-form of the differential Galois group~$H$ of the equation over the differential field $K\big(\sqrt{-1}\big)$. These facts lead us to consider two issues: determining the number of $K$-differential isomorphism classes of Picard--Vessiot extensions and describing the variation of the differential Galois group in the set of $k$-forms of~$H$. We address these two issues in the cases when~$H$ is a special linear, a special orthogonal, or a~symplectic linear algebraic group and conclude that there is no general behaviour.}

\Keywords{real Picard--Vessiot theory; linear algebraic groups; group cohomology; real forms of algebraic groups}

\Classification{12H05; 13B05; 14P05; 12D15}

\section{Introduction}
To a homogeneous linear differential equation defined over a differential field $K$ with field of constants $k$,
Picard--Vessiot theory associates a differential field extension $L$ of $K$, differentially generated over $K$
by a fundamental system of solutions of the equation, and with constant field equal to $k$, called a Picard--Vessiot extension for the given equation.
When $k$ is algebraically closed, Kolchin \cite{K} established that a Picard--Vessiot extension for the given equation exists and is unique up to $K$-differential isomorphism. The differential Galois group $\DGal(L/K)$ is defined as the group of $K$-differential automorphisms of $L$ and has the structure of a linear algebraic group defined over $k$.

For a homogeneous linear differential equation defined over a formally real differential field $K$ with real closed field of constants $k$, Crespo, Hajto and van der Put proved in \cite{CHP}, the existence and unicity up to $K$-differential isomorphism of a formally real Picard--Vessiot extension, endowed with an ordering extending the one in $K$. We note that such a linear differential equation may also have Picard--Vessiot extensions which are not formally real fields. Our result was later generalized in \cite{KP}, by using model-theoretic methods, to the case when $K$ is a differential field of characteristic 0 such that its field of constants $k$ is existentially closed in $K$ for strongly normal extensions of $K$ associated to logarithmic differential equations over $K$ on algebraic groups over~$k$. Let $L$ be a Picard--Vessiot extension of a formally real differential field $K$ with real closed field of constants $k$. Then the differential Galois group $G=\DGal(L|K)$ has the structure of a linear algebraic
group defined over $k$ (see Section~\ref{section2}).

Differential Galois theory over non-algebraically closed field of constants has been developed by several authors, see \cite{AM,An,LP}.
The inverse problem in this setting has also been considered (see \cite{BHH, D}). In particular, Dyckerhoff proved in \cite{D} that every linear algebraic group over $\R$ is a differential Galois group over the field $\R(x)$ of rational functions.

In Picard--Vessiot theory over formally real differential fields, one may find phenomena which do not arise in the context of differential fields with algebraically closed field of constants. Given a linear differential equation
\begin{gather*} \cL(Y):=Y^{(n)}+a_{n-1} Y^{(n-1)}+\dots +a_1Y'+a_0 Y=0\end{gather*}
with $a_{n-1}, \dots, a_1,a_0$ belonging to a formally real differential field $K$ with real closed field of constants $k$, one may ask the following questions which do not have a counterpart in the case when the field of constants $k$ is algebraically closed.
\begin{enumerate}\itemsep=0pt
\item[--] How many $K$-differential isomorphism classes of Picard--Vessiot extensions are there for \linebreak $\mathcal{L}(Y)=0$?
\item[--] Are the corresponding differential Galois groups $k$-isomorphic?
\end{enumerate}

Concerning the first question, as mentioned above, it is known that there is at least one Picard--Vessiot extension of $K$ for $\mathcal{L}(Y)=0$, which is a formally real field. For the second one, let us note that if $L$ and $L'$ are Picard--Vessiot extensions of a formally real differential field $K$ with real closed field of constants $k$ for the same equation, and $G=\DGal(L/K)$, $G'=\DGal(L'/K)$ are the corresponding differential Galois groups, then we have $G\times_k\overline{k}\simeq G'\times_k\overline{k}$, i.e., $G$ and $G'$ are both $k$-forms of the group $H=G\times_k\overline{k}$.

In this paper, we consider a formally real Picard--Vessiot extension $L$ of $K$ for a linear differential equation $\mathcal{L}(Y)=0$ defined over $K$ and the differential Galois group $G=\DGal(L|K)$. We give an answer to the above questions in the case in which $H:=G\times_k \overline{k}$ is a special linear, special orthogonal, or symplectic linear algebraic group. Let us note that it makes sense to start with such a Picard--Vessiot extension $L/K$ since in \cite[Proposition~3.3]{CHP} Crespo, Hajto and van der Put proved
that when the differential field $K$ is real closed, given a~connected semi-simple linear algebraic group $G$ defined over~$k$, there exists a linear differential equation defined over~$K$ and a formally real Picard--Vessiot extension $L|K$ for it such that $G=\DGal(L|K)$. The inspection of the different cases shows that there is no general pattern. The differential Galois group may be the same for all Picard--Vessiot extensions of $\mathcal{L}(Y)=0$ or range over the whole set of $k$-forms of~$H$. For a linear differential equation defined over a formally real differential field~$K$ the determination of its real differential group $G$ gives more information on the behaviour of the solutions than the determination of the complexification $H$ of $G$. For example, for $k=\R$ and $K$ a field of real functions, the determination of~$G$ will give information on the existence of oscillating functions among the solutions of $\mathcal{L}(Y)=0$ (see \cite{CH2} and the earlier topological approach in~\cite{GK}). It is then interesting to study how the real differential group $\DGal(L|K)$ varies as $L$ runs over the $K$-isomorphism classes of Picard--Vessiot extensions.

We refer the reader to \cite{CH} for the topics on differential Galois theory, when the field of constants is algebraically closed, to \cite{BCR} or \cite{P} for those on formally real fields and to \cite{Ma,Sp, St} for those on linear algebraic groups.

\section{Preliminaries}\label{section2}

For the reader's convenience, we recall the definitions of formally real field, real closed field and Picard--Vessiot extension.

\begin{Definition} A \emph{formally real field} is a field which may be given an ordering compatible with the field operations. Equivalently, a field $K$ is formally real if $-1$ is not a sum of squares in $K$.

A formally real field is \emph{real closed} if it has no nontrivial algebraic extensions which are formally real fields. Equivalently, a field $k$ is real closed if $-1$ is not a square in $k$ and $k\big(\sqrt{-1}\big)$ is algebraically closed.
\end{Definition}

We note that in the literature on real algebraic geometry ``real field'' is frequently used for ``formally real field''. A field of positive characteristic is not formally real. The field $\R$ of real numbers is the standard example of a real closed field. The field $\R(x)$ of rational functions and the field of formal Laurent series $\R(x)$ with derivation ${\rm d}/{\rm d}x$ are examples of real differential fields with real closed field of constants.

\begin{Definition} Given a linear differential equation
\begin{gather*} \cL(Y):=Y^{(n)}+a_{n-1} Y^{(n-1)}+\dots +a_1Y'+a_0 Y=0\end{gather*}
defined over a differential field $K$, with field of constants $k$, a \emph{Picard--Vessiot extension} of $K$ for
$\cL(Y)=0$ is a differential field extension $L|K$ such that
\begin{enumerate}\itemsep=0pt
\item[a)] $L$ is differentially generated over $K$ by a full set of solutions of $\cL(Y)=0$;
\item[b)] the field of constants of $L$ is $k$.
\end{enumerate}
\end{Definition}

Let us assume that $K$ is a formally real differential field with real closed field of constants $k$, $\cL(Y)=0$ a~linear differential equation defined over $K$ and $L|K$ a Picard--Vessiot extension for $\cL(Y)=0$. In this case, we note that the set $\operatorname{DHom}_K(L,L({\rm i}))$ of $K$-differential morphisms from $L$ to $L({\rm i})$ is in bijection with the set $\operatorname{DAut}_{K({\rm i})} L({\rm i})$ of $K({\rm i})$-differential automorphisms of $L({\rm i})$. We define the differential Galois group of $L|K$ as the set $\operatorname{DHom}_K(L,L({\rm i}))$ with the group structure obtained by transferring the one of $\operatorname{DAut}_{K({\rm i})} L({\rm i})$ via the above bijection. The differential Galois group has the structure of a $k$-defined linear algebraic group (see \cite[Proposition~4.1]{CHS}). We note that the proof of the existence of a formally real Picard--Vessiot extension for a linear differential equation defined over a formally real differential field with real closed field of constants given in \cite[Theorems~3.2 and~3.3 and Corollary~3.4]{CHS} is not right. The remaining proofs in~\cite{CHS} are correct.

In the sequel, $K$ will denote a formally real differential field with real closed field of constants~$k$, $\cL(Y)=0$ a
linear differential equation defined over $K$, $L|K$ a formally real Picard--Vessiot extension for $\cL(Y)=0$ and $G$ the differential Galois group of $L|K$. We
want to determine the number of Picard--Vessiot extensions of $K$ for $\cL(Y)=0$, up to $K$-differential isomorphism, and the differential Galois group of each of them.

The set of $K$-differential isomorphism classes of Picard--Vessiot extensions for $\mathcal{L}(Y)=0$ is in one-to-one correspondence
with the cohomology set $H^1\big(k,G\big(\overline{k}\big)\big)$, where $\overline{k}$ denotes the algebraic closure of $k$ and $G$ denotes the differential Galois group $\DGal(L/K)$. Indeed, we have a bijection between the set
of $K$-differential isomorphism classes of Picard--Vessiot extensions for $\mathcal{L}(Y)=0$ and the set of isomorphism classes of fiber functors $\omega\colon \langle M \rangle_{\otimes} \rightarrow \operatorname{vect}(k)$, where $\langle M \rangle_{\otimes}$ denotes the Tannakian category generated by the $K$-differential module $M$ associated to $\mathcal{L}(Y)=0$ \cite[Proposition~1]{CHP}. In turn, this set of isomorphism classes of fiber functors is in bijection with the set of isomorphism classes of $G$-torsors, by \cite[Theorem~3.2]{DM}. Finally, the set of isomorphism classes of $G$-torsors is in bijection with $H^1\big(k,G\big(\overline{k}\big)\big)$ (see, e.g., \cite[Lemma A.5.1]{PS}). If $L'$ is a~Picard--Vessiot extension for $\mathcal{L}(Y)=0$, $L({\rm i})$ and $L'({\rm i})$ are Picard--Vessiot extensions of $K({\rm i})$ for $\mathcal{L}(Y)=0$. Since the field of constants of $K({\rm i})$ is $\overline{k}$, we have an isomorphism of differential fields $f\colon L({\rm i}) \rightarrow L'({\rm i})$, by the unicity of the Picard--Vessiot extension in the case when the field of constants is algebraically closed. The group $\Gal\big(\overline{k}|k\big)$ acts on the set of isomorphisms from $L({\rm i})$ to $L'({\rm i})$ by $s(f)=s\circ f \circ s^{-1}$, for $s \in \Gal\big(\overline{k}|k\big)$, where we denote also by $s$ the automorphisms of $L({\rm i})$ and $L'({\rm i})$ induced by $s$. The $1$-cocycle $x$ corresponding to $L'$ is determined by $x(c)= f^{-1} \circ c(f)$, where $c$ is the nontrivial element in $\Gal\big(\overline{k}|k\big)$. Reciprocally, if $x$ is a 1-cocycle from $\Gal\big(\overline{k}|k\big)$ to~$G\big(\overline{k}\big)$, then the corresponding Picard--Vessiot extension corresponding to $x$ is the subfield of $L({\rm i})$ fixed by the automorphism $x(c) \circ c$, by Galois descent theory (see \cite[Chapter~III, Section~1.3]{S2}). The Galois group $\Gal\big(\overline{k}|k\big)$ acts on $G\big(\overline{k}\big)$, by an involution leaving~$G$ invariant, hence the cohomology set $H^1\big(k,G\big(\overline{k}\big)\big)$ depends on the $k$-form $G$ of $H=G\times_k \overline{k}$.

Let now $G$ denote a linear algebraic group defined over a real closed field $k$ and let $H=G\times_k \overline{k}$. The set of $k$-forms of $H$ is in one-to-one correspondence with the cohomology set $H^1\big(k,\Aut G\big(\overline{k}\big)\big)$, where $\Gal\big(\overline{k}|k\big)$ acts on $\Aut G\big(\overline{k}\big)$ by $s(f)=s\circ f \circ s^{-1}$, for $s \in \Gal\big(\overline{k}|k\big)$, $f \in \Aut G\big(\overline{k}\big)$, as usual. Let us note that the classification of the $k$-forms of $H$ is equivalent to the classification
of the real forms of the corresponding complex group.

 We consider the map
\begin{gather*} \Phi\colon \ H^1\big(k, G\big(\overline{k}\big)\big) \rightarrow H^1\big(k,\Aut G\big(\overline{k}\big)\big)\end{gather*}
induced by the morphism from $G\big(\overline{k}\big)$ to $\Aut G\big(\overline{k}\big)$ sending an element $g$ in $G\big(\overline{k}\big)$ to
conjugation by $g$. When $G$ is the differential Galois group of a~Picard--Vessiot extension $L$ of $K$ for a linear differential equation $\cL(Y)=0$, $\Phi$ sends the element in
$H^1\big(k,G\big(\overline{k}\big)\big)$ corresponding to a Picard--Vessiot extension $L'$ of $K$ for $\cL(Y)=0$ to the element in $H^1\big(k,\Aut G\big(\overline{k}\big)\big)$ corresponding to $\DGal(L'|K)$ (see \cite[Section~3, Observations~1]{CHP}).

For $c$ the nontrivial element of $\Gal\big(\overline{k}|k\big)$, we write $\overline{a}=c(a)$, for $a$ an element in $\overline{k}$. For $v=(a_1,\dots,a_n)\in \overline{k}^n$, we shall write $\overline{v}=(\overline{a_1},\dots,\overline{a_n})$ and for $M=(a_{ij})$ a matrix with entries in $\overline{k}$, $\overline{M}=(\overline{a_{ij}})$.

We shall consider a linear algebraic group $G$ defined over the real closed field $k$, such that $H=G\times_{k} \overline{k}$ is either a~special linear group $\SL(n)$, a~special
orthogonal group $\SO(n)$ or a~symplectic group~$\Sp(n)$. We will then consider the real forms of each of these groups. We refer to~\cite{S2} or~\cite{Kn} for their determination, to~\cite{B} or~\cite{T} for a~more explicit description of them.

We will determine in each case the number of $K$-differential isomorphism classes of Picard--Vessiot extensions of $K$ for $\cL(Y)=0$ and the differential Galois group for each class.

Let ${\rm i}$ denote a square root of $-1$ in $\overline{k}$. For $p$, $n$ integers with $0\leq p \leq n$, we define the $n\times n$ matrices
\begin{gather*}I_p=\left( \begin{matrix} {\rm Id}_p & 0 \\ 0 & -{\rm Id}_{n-p}
\end{matrix} \right), \qquad J_p=\left( \begin{matrix} {\rm Id}_p & 0 \\ 0 & {\rm i}\, {\rm Id}_{n-p}
\end{matrix} \right).\end{gather*}

\section[Forms of $\SL(n)$]{Forms of $\boldsymbol{\SL(n)}$}\label{section3}

The real forms of $\SL(n), n \geq 2,$ are
\begin{enumerate}\itemsep=0pt
\item[1)] $\SL(n,k)$;
\item[2)] $\SL(n/2,\h)$ if $n$ is even, where $\h$ denotes the quaternion algebra over~$k$;
\item[3)] $\SU\big(n,\overline{k},h\big)$, where $h$ is a nondegenerate hermitian form on~$\overline{k}^n$.
\end{enumerate}

For $G$ each of these real forms and $K$ a formally real differential field, with real closed field of constants $k$, we consider a linear differential equation $\cL(Y)=0$ of order $n$ defined over $K$ and a~Picard--Vessiot extension $L|K$ for $\cL(Y)=0$ such that $L$ is formally real and $\DGal(L/K) \simeq G.$

\subsection[$G=\SL(n,k)$]{$\boldsymbol{G=\SL(n,k)}$}

We have $\big|H^1\big(k,G\big(\overline{k}\big)\big)\big|=1$ \cite[Chapter~X, Section~1]{S}, hence $L|K$ is the unique Picard--Vessiot extension for $\cL(Y)=0$, up to $K$-differential isomorphism.

\subsection[$G=\SL(n/2,\h)$]{$\boldsymbol{G=\SL(n/2,\h)}$}\label{s}

Let us denote by $1$, $I$, $J$, $K$ the basis elements of $\h$. We recall that $\GL(n/2,\h)$ embeds into $\GL\big(n,\overline{k}\big)$ via the morphism $(h_{ij}) \mapsto (\mu(h_{ij}))$, where
\begin{gather*}
\mu(a+bI+cJ+dK)=\left(\begin{matrix} a+b{\rm i} &
c+d{\rm i}
\\ -c+d{\rm i} & a-b{\rm i}
\end{matrix}\right), \qquad a,b,c,d \in k.
\end{gather*}
We denote by $A_n$ the matrix $(a_{ij})_{1\leq i,j \leq
n}$ with
\begin{gather} \label{An} a_{ij}= \begin{cases} \hphantom{-}1 & \text{if} \ i \ \text{is odd and} \ j=i+1, \\ -1 &
\text{if} \ i \ \text{is even and} \ j=i-1, \\ \hphantom{-}0 & \text{in all other cases.}
\end{cases}
\end{gather}
We have $\mu(\GL(n/2,\h))=\big\{ M \in \GL\big(n,\overline{k}\big)\colon M=A_n \overline{M} A_n^{-1} \big\}$ and $\mu(\SL(n/2,\h))= \SL\big(n,\overline{k}\big) \cap \mu(\GL(n/2,\h))$.

By \cite[Chapter~VII, Section~29, Corollary~29.4]{Inv}, we have $H^1\big(k,G\big(\overline{k}\big)\big) \simeq k^* /{\rm Nrd}\big(\h^{n/2}\big)$. Since the norm of a quaternion is always positive, we obtain $\big|H^1\big(k,G\big(\overline{k}\big)\big)\big|=2$ (see also \cite[Chapter~III, Section~1.4]{S2}). We have then two Picard--Vessiot
extensions for $\cL(Y)=0$, up to $K$-differential isomorphism. A nontrivial 1-cocycle $x$ of $\Gal\big(\overline{k}|k\big)$ in
$\SL\big(n,\overline{k}\big)$ is given by $x(c)= \zeta\,{\rm Id}$, for $\zeta$ a~primitive $n$-th root of unity. A
$K$-differential automorphism of $L({\rm i})$ corresponding to~$x$ is
given by the matrix $A:=\zeta^{-1/2} \, {\rm Id}$ on the vector space of
solutions, since $A$ satisfies $A^{-1} c(A)=x(c)$ (see \cite[Chapter~X, Section~2, Proposition~4]{S}).
Conjugation by $\zeta^{-1/2} \, {\rm Id}$ leaves the group $\SL(n/2,\h)$ stable. We obtain that the Picard--Vessiot extensions
for $\cL(Y)=0$ in both $K$-differential isomorphy classes have $\SL(n/2,\h)$ as differential Galois group.

\subsection[$G=\SU\big(n,\overline{k},h\big)$]{$\boldsymbol{G=\SU\big(n,\overline{k},h\big)}$}\label{SU}

It is known that if $h$ is a nondegenerate hermitian form on $\overline{k}^n$, then $h$ is equivalent to a hermitian form with matrix $I_p$, for some integer $p$ with $0\leq p \leq n$, called the index of $h$ and that two nondegenerate hermitian forms on $\overline{k}^n$ are equivalent if and only if their indices coincide (see \cite[Chapter~VII, Section~29, Example~29.19]{Inv} and \cite[Section~3.3]{BP}).

We fix $G=\big\{ M \in \SL\big(n,\overline{k}\big)\colon \overline{M}^t I_p M = I_p \big\}$ and consider the action of $\Gal\big(\overline{k}|k\big)$ on $\SL\big(n,\overline{k}\big)$ given by $c(M)=I_p \big(\overline{M}^t\big)^{-1} I_p$. We shall prove
\begin{gather}
\big|H^1\big(k,G\big(\overline{k}\big)\big)\big|= \begin{cases} \left[ \dfrac n 2 \right]+1, & \text{when $n$ is odd or $p$ is even},
 \vspace{2mm}\\
 \dfrac n 2, & \text{when $n$ is even and $p$ is odd}.\end{cases}
\label{h1}
\end{gather}

To this end we shall determine a maximal set of pairwise nonequivalent $1$-cocycles from $\Gal\big(\overline{k}|k\big)$ in $\SL\big(n,\overline{k}\big)$.
For such a cocycle $x$, we may assume that the image of ${\rm Id} \in \Gal\big(\overline{k}|k\big)$ is the identity matrix and then $x$ is determined by the image $B \in \SL\big(n,\overline{k}\big)$ of the unique non trivial element $c$ in $\Gal\big(\overline{k}|k\big)$. By the $1$-cocycle condition, $B$ must satisfy
$Bc(B)= {\rm Id}$. We denote by $x_q$ the cocycle given by $c \mapsto
B_q$, where $ B_q=I_q I_p$, with $q$ an integer of the same
parity as $p$ and $0\leq q \leq n$. Let us see that every 1-cocycle $x$ from
$\Gal\big(\overline{k}|k\big)$ in $\SL\big(n,\overline{k}\big)$ is equivalent to
some $x_q$. As said above, such a cocycle $x$ is determined by $x(c)=B$
satisfying $Bc(B)={\rm Id}$, i.e., $BI_p\big(\overline{B}^t\big)^{-1}I_p={\rm Id}$, equivalently
$BI_p=I_p\overline{B}^t=\big(\overline{BI_p}\big)^t$, so $BI_p$ is an
hermitian matrix, hence there exists an invertible matrix $M$ such
that $M^{-1} BI_p\big(\overline{M}^t\big)^{-1}=I_r$, for some integer $r$, $0\leq r \leq n$.
Equivalently
\begin{gather}\label{eq} M^{-1} B c(M)=I_r I_p. \end{gather}
Let \looseness=-1 us note that, taking determinants in (\ref{eq}), we obtain $\det M\, \overline{\det M}=1$, hence there exists $\zeta \in \overline{k}$ such that
$\det(\zeta M)=1$ and $\zeta M$ satisfies~(\ref{eq}). We have then that there exists a matrix $M \in \SL\big(n,\overline{k}\big)$
satisfying (\ref{eq}), which means that $x$ is equivalent to the $1$-cocycle $x_r$ determined by $c \mapsto I_rI_p$. Since $B \in \SL\big(n,\overline{k}\big)$, we have $\det(I_r
I_p)=1$, so $r$ is an integer of the same parity as~$p$.

Let us see now that the 1-cocycles $x_q$ are pairwise nonequivalent.
We have $x_q \sim x_{q'} \Leftrightarrow \exists\, M \in \SL\big(n,\overline{k}\big)$ such that $B_{q'}=M^{-1} B_q c(M)$. This
equality is equivalent to $MI_{q'}\overline{M}^t=I_q$ which
implies $q=q'$, so the 1-cocycles $x_q$ are pairwise nonequivalent. We have then $H^1\big(k,G\big(\overline{k}\big)\big)=\{[x_q] \colon 0\leq q \leq n, \, q\equiv p \pmod{2} \}$.
Then $\big|H^1\big(k,G\big(\overline{k}\big)\big)\big|=|\{q \in \Z \colon 0\leq q \leq n, q\equiv p \pmod{2} \}|$ and we obtain the values in (\ref{h1}).

We have $Z\big(\SL\big(n,\overline{k}\big)\big)=\mu_n\big(\overline{k}\big)$. We want to
determine the image of $[x_q]$ under the map
$\Phi\colon H^1\big(k,G\big(\overline{k}\big)\big) \rightarrow H^1\big(k,\Aut G\big(\overline{k}\big)\big)$. The 1-cocycle $x_q$ corresponds to a~Picard--Vessiot extension $L_q$ of $K$ for $\cL(Y)=0$ such that there is a differential isomorphism $f_q$ from $L({\rm i})$ to $L_q({\rm i})$ satisfying $x_q=f_q^{-1} c(f_q)$. The differential isomorphism $f_q$ is determined by the matrix $D_q$ giving the images of a vector space of solutions. The isomorphism $f_q$ satisfies $x_q=f_q^{-1} c(f_q)$ if and only if the matrix $D_q$ satisfies $B_q=D_q^{-1} \overline{D_q}$. We may take $D_q:=J_qJ_p$. Since conjugation by $D_q$ leaves the group $G$ invariant, we obtain that all Picard--Vessiot extensions of $K$ for $\cL(Y)=0$ have the same differential Galois group~$G$.

Gathering the results in this section we may state the following theorem.

\begin{Theorem} Let $K$ be a formally real differential field with real closed field of constants $k$, $\cL(Y)=0$ a
linear differential equation defined over $K$, $L|K$ a formally real Picard--Vessiot extension for $\cL(Y)=0$ and $G$ the differential Galois group of $L|K$. We assume that $G$ is a real form of $\SL(n)$.
\begin{enumerate}\itemsep=0pt
\item[$(1)$] If $G=\SL(n,k)$, $L|K$ is the unique Picard--Vessiot extension for the equation $\cL(Y)=0$, up to $K$-differential isomorphism.
\item[$(2)$] If $G=\SL(n/2,\h)$, there are two Picard--Vessiot extensions for the equation $\cL(Y)=0$, up to $K$-differential isomorphism, and both of them have differential Galois group $G$.
\item[$(3)$] If $G=\SU\big(n,\overline{k},h\big)$, there are $[n/2]+1$ $($resp.~$[n/2])$ Picard--Vessiot extensions for the equation $\cL(Y)=0$, if $n$ is odd or $p$ is even $($resp.\ if $n$ is even and $p$ is odd$)$, up to $K$-differential isomorphism, and all of them have differential Galois group~$G$.
\end{enumerate}
\end{Theorem}

\section[Forms of $\SO(n)$]{Forms of $\boldsymbol{\SO(n)}$}\label{section4}

The real forms of $\SO(n)$, with $n$ odd, are the groups $\SO(n,k,Q)$, where $Q$ is a nondegenerate quadratic form on $k^{n}$. When $G$ is one of these forms, we proved in \cite[Section~3, Examples~1 and~3]{CHP} that the map $\Phi\colon H^1\big(k, G\big(\overline{k}\big)\big) \rightarrow H^1\big(k,\Aut G\big(\overline{k}\big)\big)$ is a bijection. We consider now the case when $n$ is even.

The real forms of $\SO(n)$, with $n$ even, are
\begin{enumerate}\itemsep=0pt
\item[1)] $\SO(n,k,Q)$, where $Q$ is a nondegenerate quadratic form on
$k^{n}$;
\item[2)] $\SU(n/2,\h,h)$, where $h$ is a nondegenerate anti-hermitian form on $\h^{n/2}$ (with respect to the involution $\sigma$ of $\h$ defined by
$a+bI+cJ+dK \mapsto a-bI-cJ-dK$).
\end{enumerate}

For $G$ each of these real forms and $K$ a formally real differential field, with real closed field of constants $k$, we consider a linear differential equation $\cL(Y)=0$ of order $n$ defined over $K$ and a~Picard--Vessiot extension $L|K$ for $\cL(Y)=0$ such that $L$ is formally real and $\DGal(L/K) \simeq G.$

\subsection[$G=\SO(n,k,Q)$]{$\boldsymbol{G=\SO(n,k,Q)}$}

The quadratic form $Q$ is equivalent to a quadratic form with matrix $I_p$, for some integer $p$ with $0\leq p \leq n$, which determines the equivalence class of $Q$.

The cohomology set $H^1\big(k,G\big(\overline{k}\big)\big)$ is in one-to-one correspondence with the set of equivalence classes of quadratic forms on $k^n$ of rank $n$ and index of the same parity as $p$ (see \cite[Chapter~VII, Section~29, formula~(29.29)]{Inv}). We
have then
\begin{gather*} \big|H^1\big(k,G\big(\overline{k}\big)\big)\big|= \begin{cases} \dfrac n 2 +1, & \text{when $p$ is even},\vspace{2mm}\\
\dfrac n 2, & \text{when $p$ is odd}.\end{cases} \end{gather*}

The cocycles $x_q$ defined by $c \mapsto B_q$, where $ B_q=I_q I_p$, with $q$ an integer of the same parity as~$p$ and $0\leq q
\leq n$, form a complete system of representatives of the cohomology set $H^1\big(k,G\big(\overline{k}\big)\big)$. The 1-cocycle $x_q$ corresponds to a~Picard--Vessiot extension $L_q$ of $K$ for $\cL(Y)=0$ such that there is a differential isomorphism $f_q$ from $L({\rm i})$ to $L_q({\rm i})$ satisfying $x_q=f_q^{-1}
c(f_q)$. The differential isomorphism $f_q$ is determined by the matrix $D_q$ giving the images of a vector space of solutions. The isomorphism $f_q$ satisfies $x_q=f_q^{-1}
c(f_q)$ if and only if the matrix $D_q$ satisfies $B_q=D_q^{-1} \overline{D_q}$. We may take $D_q:=J_qJ_p$. If the matrix $M$ satisfies $M^t I_p M=I_p$, the conjugate matrix $N:=D_qMD_q^{-1}$ satisfies $N^t I_q N=I_q$, hence the Picard--Vessiot extension corresponding to the 1-cocycle $x_q$ has differential Galois group $\SO(n,k,Q_q)$, where $Q_q$ denotes the quadratic form with index $q$. Let us note that $\SO(n,k,Q_q)=\SO(n,k,Q_{n-q})$, hence the Picard--Vessiot extension corresponding to $x_q$ and $x_{n-q}$, $0\leq q \leq (n/2)-1$, have the same differential Galois group.

\subsection[$G=\SU(n/2,\h,h)$, $h$ anti-hermitian]{$\boldsymbol{G=\SU(n/2,\h,h)}$, $\boldsymbol{h}$ anti-hermitian}

We have $\U(n/2,\h,h)=\big\{ M \in \GL(n/2,\h)\colon \sigma(M)^t [h] M = [h]\big\}$, for $[h]$ the matrix of the anti-hermitian form $h$, in some basis of $\h^{n/2}$. The group $\U(n/2,\h,h)$ is the group of automorphisms of the anti-hermitian vector space $\big(\h^{n/2},h\big)$. The set of equivalence classes of nondegenerate anti-hermitian forms over $\h^{n/2}$ is in one-to-one correspondence with the cohomology set $H^1\big(k,\U(n/2,\h,h)\big(\overline{k}\big)\big)$. Up to equivalence, there is one single anti-hermitian form on $\h^n$, hence $H^1\big(k,\U(n/2,\h,h)\big(\overline{k}\big)\big)=1$.

We consider the exact sequence
\begin{gather*} 1 \rightarrow \SU \rightarrow \U \rightarrow \mu_2 \rightarrow 1.\end{gather*}
Since the reduced norm of a quaternion is always positive, the reduced norm $\U(n/2,\h,h)(k) \rightarrow \mu_2(k)$ is the trivial map. We obtain then for the cohomology sets the exact sequence
\begin{gather*} 1 \rightarrow \mu_2 \rightarrow H^1\big(k,\SU(n/2,\h,h)\big(\overline{k}\big)\big) \rightarrow H^1\big(k,\U(n/2,\h,h)\big(\overline{k}\big)\big).\end{gather*}
Therefore $\big|H^1\big(k,\SU(n/2,\h,h)\big(\overline{k}\big)\big)\big|=2$. We obtain then that there are two Picard--Vessiot extensions for $\cL(Y)=0$, up to $K$-differential isomorphism. We denote by~$L'$ the non formally real one. We may check that $\mu(G)$ is the intersection of (a conjugate form of) $\SO\big(n,\overline{k}\big)$ with $\mu(\GL(n/2,\h))$. A nontrivial 1-cocycle of $\Gal\big(\overline{k}|k\big)$ in $G\big(\overline{k}\big)$ is given by $c\mapsto A_n$, for $A_n$ the matrix defined by~(\ref{An}). The matrix $B=(b_{ij})$ defined by
\begin{gather*}b_{ij}= \begin{cases} \hphantom{-}\dfrac{1}{\sqrt{2}} & \text{if } \ i=j \ \text{or} \ i \ \text{is even and} \ j=i-1, \vspace{1mm}\\
-\dfrac{1}{\sqrt{2}} &
\text{if} \ i \ \text{is odd and} \ j=i+1, \vspace{1mm}\\ \hphantom{-} 0 & \text{in all other cases}.
\end{cases}
\end{gather*}
satisfies $B^{-1}c(B)=A_n$, hence the cohomology class $[x]$ corresponds to the isomorphism class of the differential isomorphism $f$ from $L({\rm i})$ to $L'({\rm i})$ with matrix $B$ on the vector space of solutions. Since conjugation by $B$ leaves $G$ invariant, we obtain that both Picard--Vessiot extensions have the same differential Galois group.

Gathering the results in this section we may state the following theorem. For completeness, we include the case $n$ odd.

\begin{Theorem} Let $K$ be a formally real differential field with real closed field of constants $k$, $\cL(Y)=0$ a~linear differential equation defined over $K$, $L|K$ a formally real Picard--Vessiot extension for $\cL(Y)=0$ and $G$ the differential Galois group of $L|K$. We assume that $G$ is a real form of $\SO(n)$.

If $n$ is odd, there are $(n+1)/2$ Picard--Vessiot extensions for the equation $\cL(Y)=0$, up to $K$-differential isomorphism, and their differential Galois groups range over the whole set of real forms of $\SO(n)$.

 Assume $n$ even.
\begin{enumerate}\itemsep=0pt
\item[$(1)$] If $G=\SO(n,k,Q_p)$, where $Q_p$ is a nondegenerate quadratic form on $k^n$, of index $p$, there are $(n/2) +1$ $($resp.~$n/2)$ Picard--Vessiot extensions for the equation $\cL(Y)=0$, up to $K$-differential isomorphism, when $p$ is even $($resp.\ when $p$ is odd$)$ and their differential Galois groups range over the whole set of groups $G=\SO(n,k,Q_q)$, with $Q_q$ a nondegenerate quadratic form on~$k^n$, of index~$q$, $0\leq q \leq n/2$ and $q$ of the same parity as $p$.
\item[$(2)$] If $G=\SU(n/2,\h,h)$, where $h$ is a nondegenerate anti-hermitian form on~$\h^n$, there are two Picard--Vessiot extensions for the equation $\cL(Y)=0$, up to $K$-differential isomorphism, and they have both differential Galois group~$G$.
\end{enumerate}
\end{Theorem}

\section[Forms of $\Sp(2n)$]{Forms of $\boldsymbol{\Sp(2n)}$}\label{section5}

The real forms of $\Sp(2n)$ are
\begin{enumerate}\itemsep=0pt
\item[1)] $\Sp(2n,k)$;
\item[2)] $\SU(n,\h,h)$, where $h$ is a nondegenerate hermitian form on $\h^n$ (with respect to the involution $\sigma$ of $\h$ defined by $a+bI+cJ+dK \mapsto a-bI-cJ-dK$).
\end{enumerate}

For $G$ each of these real forms and $K$ a formally real differential field, with real closed field of constants $k$, we consider a linear differential equation $\cL(Y)=0$ of order $n$ defined over $K$ and a~Picard--Vessiot extension $L|K$ for $\cL(Y)=0$ such that $L$ is formally real and $\DGal(L/K) \simeq G.$

\subsection[$G=\Sp(2n,k)$]{$\boldsymbol{G=\Sp(2n,k)}$}

We have $\big|H^1\big(k,G\big(\overline{k}\big)\big)\big|=1$ \cite[Chapter~X, Section~2, Corollary~2]{S}, hence $L|K$ is the unique Picard--Vessiot extension for
$\cL(Y)=0$, up to $K$-differential isomorphism.

\subsection[$G=\SU(n,\h,h)$, $h$ hermitian]{$\boldsymbol{G=\SU(n,\h,h)}$, $\boldsymbol{h}$ hermitian}

If $h$ is a nondegenerate hermitian form on $\h^n$, then $h$ is equivalent to a hermitian form with matrix $I_p$, with $p\geq n-p$~\cite[Section~7.5.3]{B}. The number of equivalence classes of nondegenerate hermitian forms over $\h^n$ is then $\big[ \frac n 2 \big] +1$. We fix
\begin{gather*}G=\big\{ M \in \GL(n,\h) \colon \sigma(M)^tI_pM = I_p\big\}.\end{gather*}
The group $G$ is the group of automorphisms of the
hermitian vector space $\big(\h^n,h\big)$. Hence the set
of equivalence classes of nondegenerate hermitian forms over
$\h^n$ is in one-to-one correspondence with the cohomology set
$H^1\big(k,G\big(\overline{k}\big)\big)$ (see \cite[Chapter~X, Section~2, Proposition~4]{S} or \cite[Section~2.6, Lemma~3]{Kn}). Since the set
of $K$-differential isomorphism classes of Picard--Vessiot extensions for $\cL(Y)=0$
is also in one-to-one correspondence with the cohomology set
$H^1\big(k,G\big(\overline{k}\big)\big)$, we have that the number of $K$-differential isomorphism classes of Picard--Vessiot extensions for $\cL(Y)=0$ is equal to the number of equivalence classes of nondegenerate hermitian forms over $\h^n$. Let us note that for $M,N \in M_n(\h)$, we have $\sigma(MN)^t=\sigma(N)^t \sigma(M)^t$. We determine now the image of $G$ under the isomorphism $\mu\otimes_k\overline{k}\colon \GL\big(n,\h\otimes_k\overline{k}\big) \rightarrow \GL\big(2n,\overline{k}\big)$. For $M \in \GL\big(n,\h\otimes_k \overline{k}\big)$, we have $\big(\mu\otimes_k \overline{k}\big)\big(\sigma\big(M^t\big)\big)=A_{2n} \big(\mu\otimes_k \overline{k}\big)(M)^t
A_{2n}^{-1}$, for $A_{2n}$ the matrix $(a_{ij})_{1\leq i,j \leq 2n}$ defined by
\begin{gather*}a_{ij} = \begin{cases} \hphantom{-}1 & \text{if} \ i \ \text{is odd and} \ j=i+1, \\
-1 & \text{if} \ i \ \text{is even and} \ j=i-1, \\
\hphantom{-}0 & \text{in all other cases}. \end{cases}\end{gather*}
Hence $\sigma(M)^tI_pM = I_p$ implies $\big(\mu\otimes_k
\overline{k}\big)(M)^t A_{2n}^{-1} I_{2p} \big(\mu\otimes_k
\overline{k}\big)(M)= A_{2n}^{-1} I_{2p}$. Now the group
\begin{gather*}\overline{G}:= \big\{ N \in \GL\big(2n,\overline{k}\big)\colon N^t A_{2n}^{-1} I_{2p} N= A_{2n}^{-1}
I_{2p} \big\} \end{gather*}
is a conjugate form of $\Sp\big(2n,\overline{k}\big)$ and
$\mu(G)=\big\{ N \in \overline{G}\colon N=A_{2n} \overline{N} A_{2n}^{-1} \big\}$. A complete set of nonequivalent 1-cocycles of $\Gal\big(\overline{k}|k\big)$ in $\mu(G)$ is given by
\begin{align*} x_q\colon \ \Gal\big(\overline{k}|k\big)& \rightarrow \mu(G), \\
 c & \mapsto B_q:=I_{2q} I_{2p} ,\end{align*}
with $q$ an integer, $0\leq q \leq n$, $q\geq n-q$. The 1-cocycle $x_q$ corresponds to a~Picard--Vessiot extension $L_q$ of $K$ for $\cL(Y)=0$ such that there is a differential isomorphism $f_q$ from~$L({\rm i})$ to~$L_q({\rm i})$ satisfying $x_q=f_q^{-1} c(f_q)$. The differential isomorphism $f_q$ is determined by the matrix~$D_{2q}$ giving the images of a vector space of solutions. The isomorphism $f_q$ satisfies $x_q=f_q^{-1} c(f_q)$ if and only if the matrix $D_{2q}$ satisfies $B_q=D_{2q}^{-1} \overline{D_{2q}}$. We may take $D_{2q}:=J_{2q}J_{2p}$. If $N$ is a~matrix belonging to~$\mu(G)$, it satisfies $N^t A_{2n}^{-1} I_{2p} N= A_{2n}^{-1} I_{2p}$. Then the conjugate matrix~$P$ of~$N$ by $D_{2q}$, $P:=D_{2q}ND_{2q}^{-1}$, satisfies $P^t A_{2n}^{-1} I_{2q} P= A_{2n}^{-1} I_{2q}$, hence the Picard--Vessiot extension corresponding to the 1-cocycle $x_q$ has differential Galois group $\SU(n,\h,h_q)$, where $h_q$ denotes the hermitian form with index~$q$.

Gathering the results in this section we may state the following theorem.

\begin{Theorem} Let $K$ be a formally real differential field with real closed field of constants $k$, $\cL(Y)=0$ a
linear differential equation defined over $K$, $L|K$ a formally real Picard--Vessiot extension for $\cL(Y)=0$ and $G$ the differential Galois group of $L|K$. We assume that $G$ is a real form of $\Sp(2n)$.
\begin{enumerate}\itemsep=0pt
\item[$(1)$] If $G=\Sp(2n,k)$, $L|K$ is the unique Picard--Vessiot extension for the equation $\cL(Y)=0$.
\item[$(2)$] If $G=\SU(n,\h,h_p)$, where $h_p$ is a nondegenerate hermitian form on $\h^n$, of index $p$, $0\leq p \leq n$, $p\geq n-p$, there are $[n/2]+1$ Picard--Vessiot extensions for the equation $\cL(Y)=0$, up to $K$-differential isomorphism, and their differential Galois groups range over the whole set of groups $G=\SU(n,\h,h_q)$, with $h_q$ a nondegenerate hermitian form on $\h^n$, of index $q$, $0\leq q \leq n$, $q\geq n-q$.
\end{enumerate}
\end{Theorem}

\section{Conclusions}\label{section6}

In the preceding we have seen cases in which a linear differential equation $\cL(Y)=0$ defined over a formally real differential field $K$ has Picard--Vessiot extensions which are not formally real. The occurrence of these extensions depends on the real form of the differential Galois group of $\cL(Y)=0$. When the number of $K$-differential isomorphisms of Picard--Vessiot extensions of $\cL(Y)=0$ is bigger than~1, we find several situations concerning the differential Galois group, either it is the same for all Picard--Vessiot extensions or it ranges over a subset or the whole set of real forms of the group~$G$ of the formally real Picard--Vessiot extension. It would be interesting to know if, in the case when~$K$ is a field of real functions, the solutions of such an equation in a non formally real differential field and the variation of the differential Galois group have some physical interpretation. Some inspiring examples in Hamiltonian mechanics are presented in~\cite{Au}.

\subsection*{Acknowledgements} Both authors acknowledge support of grant MTM2015-66716-P (MINECO/FEDER, UE). The authors thank the anonymous referees for their valuable remarks and suggestions.

\pdfbookmark[1]{References}{ref}
\LastPageEnding

\end{document}